\documentclass[preprint,12pt]{elsarticle}

\usepackage{graphics}
\usepackage{graphicx}
\usepackage{amsmath}
\usepackage{amsfonts}
\usepackage{amssymb}
\usepackage{comment}
\usepackage{rotating}
\usepackage{color}

\newcommand{\R}{\mathbb{R}}

\newtheorem{theorem}{Theorem}
\newtheorem{lemma}{Lemma}
\newtheorem{definition}{Definition}
\newtheorem{example}{Example}
\newtheorem{remark}{Remark}
\newtheorem{corollary}{Corollary}

\journal{Linear Algebra and its Applications}

\begin{document}

\begin{frontmatter}

\title{On a Construction of Integrally Invertible Graphs and their Spectral Properties}

\author[auth1]{So\v{n}a Pavl\'{\i}kov\'a}
\author[auth2]{Daniel \v{S}ev\v{c}ovi\v{c}}
\address[auth1]{Institute of Information Engineering, Automation, and Mathematics, 
FCFT, Slovak Technical University, 812 37 Bratislava, Slovakia, {\tt sona.pavlikova@stuba.sk}}
\address[auth2]{Department of Applied Mathematics and Statistics, FMFI, Comenius University, 842 48 Bratislava, Slovakia, {\tt sevcovic@fmph.uniba.sk}}

\begin{abstract}

Godsil (1985) defined a graph to be invertible if it has a non-singular adjacency matrix whose inverse is diagonally similar to a nonnegative integral matrix; the graph defined by the last matrix is then the inverse of the original graph. In this paper we call such graphs positively invertible and introduce a new concept of a negatively invertible graph by replacing the adjective `nonnegative' by `nonpositive in Godsil's definition; the graph defined by the negative of the resulting matrix is then the negative inverse of the original graph. We propose new constructions of integrally invertible graphs (those with non-singular adjacency matrix whose inverse is integral) based on an operation of `bridging' a pair of integrally invertible graphs over subsets of their vertices, with sufficient conditions for their positive and negative invertibility. We also analyze spectral properties of graphs arising from bridging and derive lower bounds for their least positive eigenvalue. As an illustration we present a census of graphs with a unique 1-factor on $m\le 6$ vertices and determine their positive and negative invertibility.

\end{abstract}


\end{frontmatter}

\section{Introduction}

A number of ways of introducing inverses of graphs have been proposed, all based on inverting adjacency matrices. For a graph with a non-singular adjacency matrix a first thought might be to hope that the inverse matrix defines a graph again. It turns out, however, that this happens to be the case only for unions of isolated edges \cite{Har}. A successful approach was initiated by Godsil \cite{Godsil1985} who defined a graph to be invertible if the inverse of its (non-singular) adjacency matrix is diagonally similar (c.f. \cite{Zas}) to a nonnegative integral matrix representing the adjacency matrix of the inverse graph in which positive labels determine edge multiplicities. This way of introducing invertibility has the appealing property that inverting an inverse graph gives back the original graph. For a survey of results and other approaches to graph inverses we recommend \cite{McMc}.

%%%%%%%

Inverse graphs are of interest in estimating the least positive eigenvalue in families of graphs, a task for which there appears to be lack of suitable bounds. However, if the graphs are invertible, one can apply one of the (many) known upper bounds on largest eigenvalues of the inverse graphs instead (cf. \cite{Pavlikova1990, Pavlikova2015}). Properties of spectra of inverse graphs can also be used to estimate the difference between the minimal positive and maximal negative eigenvalue (the so-called HOMO-LUMO gap) for structural models of chemical molecules, as it was done e.g. for graphene in \cite{YeKlein}.

%%%%%%%%

Godsil's ideas have been further developed in several ways. Akbari and  Kirkland \cite{KirklandAkb2007} and Bapat and Ghorbani \cite{Bapat} studied inverses of edge-labeled graphs with labels in a ring, Ye \emph{et al.} \cite{Ye} considered connections of graph inverses with median eigenvalues, and Pavl\'{\i}kov\'a \cite{Pavlikova2015} developed constructive methods for generating invertible graphs by edge overlapping. A large number of related results, including a unifying approach to inverting graphs, were proposed in a recent survey paper by McLeman and McNicholas \cite{McMc}, with emphasis on inverses of bipartite graphs and diagonal similarity to nonnegative matrices.

%%%%%%%%%

Less attention has been given to the study of invertibility of non-bipartite graphs and their spectral properties which is the goal of this paper. After introducing basic concepts, in Section 2 we present an example of a non-bipartite graph representing an important chemical molecule of fulvene. Its adjacency matrix has the remarkable additional property that its inverse is integral and diagonally similar to a nonpositive rather than a nonnegative matrix. This motivated us to introduce negative invertibility as a natural counterpart of Godsil \cite{Godsil1985} concept: A negatively invertible graph is one with a non-singular adjacency matrix whose inverse is diagonally similar to a nonpositive matrix. The negative of this matrix is then the adjacency matrix of the inverse graph. Positively and negatively invertible graphs are subfamilies of integrally invertible graphs, whose adjacency matrices have an integral inverse. The corresponding inverse graphs, however, would have to be interpreted as labeled graphs with both positive and negative (integral) edge labels.

%%%%%%% 

Results of the paper are organized as follows. In Section 3 we develop constructions of new integrally invertible graphs from old ones by `bridging' two such graphs over subsets of their vertices. This yields a wide range of new families of integrally invertible graphs. We derive sufficient conditions for their positive and negative invertibility. In contrast to purely graph-theoretical approach we use methods of matrix analysis and in particular results on inverting block matrices such as the Schur complement theorem and the Woodbury and Morrison-Sherman formulae. Using this approach enables us to derive useful bounds on the spectra of graphs arising from bridging construction in Section 4. We then illustrate our results in Section 5 on a recursively defined family of fulvene-like graphs. In the final Section 6 we discuss arbitrariness in the bridging construction and give a census of all invertible graphs on at most 6 vertices with a unique 1-factor.

\section{Invertible graphs}\label{sec-intro}

In this section we recall a classical concept of an invertible graph due to Godsil \cite{Godsil1985}. Let $G$ be an undirected graph, possibly with multiple edges, and with a (symmetric) adjacency matrix $A_G$. Conversely, if $A$ is a nonnegative integral symmetric matrix, we will use the symbol $G_A$ to denote the graph with the adjacency matrix $A$.

The spectrum $\sigma(G)$ of $G$ consists of eigenvalues (i.e., including multiplicities) of $A_G$ (cf. \cite{Cvetkovic1978, Cvetkovic1988}). If the spectrum does not contain zero then the adjacency matrix $A$ is invertible. We begin with a definition of an integrally invertible graph.

\begin{definition}\label{def-integerinv}
A graph $G=G_A$ is said to be integrally invertible if the inverse $A^{-1}$ of its adjacency matrix exists and is integral.
\end{definition}

%%%%%%%

It follows that a graph $G_A$ is integrally invertible if and only if $det(A)=\pm1$ (cf. \cite{KirklandAkb2007}). Note that, in such a case the inverse matrix $A^{-1}$ need not represent a graph as it may contain negative entries.

Following the idea due to Godsil, the concept of the inverse graph $G_A^{-1}$ is based on the inverse matrix $A^{-1}$ for which we require signability to a nonnegative or nonpositive matrix. We say that a symmetric matrix $B$ is positively (negatively) {signable} if there exists a diagonal $\pm 1$ matrix $D$ such that $D B D$ is nonnegative (nonpositive). We also say that $D$ is a signature matrix.

\begin{definition}\label{def:inv}
A  graph $G_A$ is called positively (negatively) invertible  if $A^{-1}$ exists and is signable to a nonnegative (nonpositive) integral matrix. If $D$ is the corresponding signature matrix, the positive (negative) inverse graph $H=G_A^{-1}$ is defined by the adjacency matrix $A_H=D A^{-1} D$ ($A_H= - DA^{-1} D$).
\end{definition}

The concept of positive invertibility coincides with the original notion of invertibility introduced by Godsil \cite{Godsil1985}. Definition~\ref{def:inv} extends Godsil's original  concept to a larger class of integrally invertible graphs with inverses of adjacency matrices signable to nonpositive matrices. 

Notice that for a diagonal matrix $D^A$ containing  $\pm1$ elements only, we have $D^A D^A = I$. It means that $(D^A)^{-1} = D^A$.

%%%%%%%%%%%%%%%%%%%%%

\begin{remark}
The idea behind the definition of an inverse graph consists of the following useful property. If $G$ is a positively (negatively) invertible graph then $G^{-1}$ is again a positively (negatively) invertible graph and $G=(G^{-1})^{-1}$. As far as the spectral properties are concerned, we have
\[
 \sigma(G^{-1}) = 1/\sigma(G) = \{ 1/\lambda, \  \lambda\in\sigma(G) \},
\]
for any positively invertible graph $G$. On the other hand, if $G$ is negatively invertible then
 \[
 \sigma(G^{-1}) = -1/\sigma(G) = \{ - 1/\lambda, \  \lambda\in\sigma(G) \}.
\]

\end{remark}

\begin{figure}
\begin{center}
\includegraphics[width=3.5truecm]{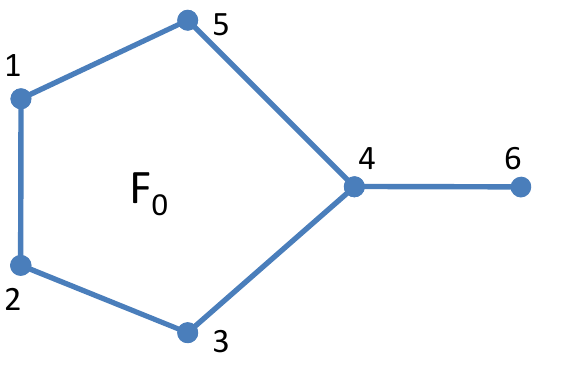}
\ \ 
\includegraphics[width=3.5truecm]{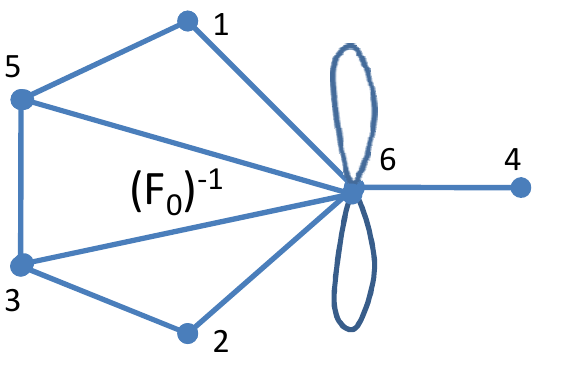}
\ \ 
\includegraphics[width=2.5truecm]{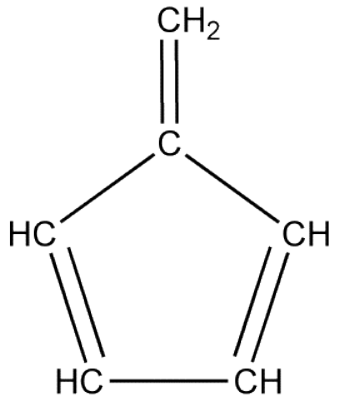}
\end{center}
\caption{
An example of a negatively invertible non-bipartite graph $F_0$ (left) and its inverse graph $(F_0)^{-1}$ (middle) of the fulvene chemical organic molecule (right).}
\label{fig-fulvene}
\end{figure}

Fig.~\ref{fig-fulvene} (left) shows the graph $F_0$ on 6 vertices representing the organic molecule of the fulvene hydrocarbon (5-methylidenecyclopenta-1,3-diene) (right). The graph $F_0$ is negatively (but not positively) invertible with the inverse graph $(F_0)^{-1}$ depicted in Fig.~\ref{fig-fulvene} (middle). The spectrum consists of the following eigenvalues:
\[
 \sigma(F_0)= \{ -1.8608, -q, -0.2541, 1/q, 1, 2.1149 \},
\]
where $q=(\sqrt{5}+1)/2$ is the golden ratio with the least positive eigenvalue $\lambda_1^+(F_0)=1/q$.  The inverse adjacency matrix $A_{F_0}^{-1}$ is signable to a nonpositive integral matrix by the signature matrix $D^{A_{F_0}} = diag(-1,-1,1,1,1-1)$.

\subsection{Bipartite graphs and their invertibility}

A graph $G_B$ is called bipartite if the set of vertices can be partitioned into two disjoint subsets such that no two vertices within the same subset are adjacent. The adjacency matrix $B$ of a bipartite graph $G_B$ can be given in a block form:
\[
 B = 
\left( 
\begin{array}{cc}
0 & K  \\
K^T & 0\\
\end{array}
\right),
\]
where $K$ is a matrix with nonnegative integer entries. Clearly, the adjacency matrix $B$ of a bipartite graph $G_B$ is invertible if and only if the number of its vertices is even and the matrix $K$ is invertible. If we consider the labeled graph corresponding to the adjacency matrix $B^{-1}$ and the product of edge labels on every cycle in this graph is positive, then $B^{-1}$ is signable to a nonnegative matrix and so $G_B$ is a positively invertible graph (see \cite{KirklandTif2009}).

Recall that a 1-factor (a perfect matching) of a graph is a spanning 1-regular subgraph with all vertices of degree 1. If $G$ is a bipartite graph with a 1-factor $M$  such that the graph $G/M$ obtained from $G$ by contracting edges of $M$ is bipartite then $G$ is a positively invertible graph (c.f. \cite{Godsil1985, Pavlikova1990}). 

Bipartiteness and invertibility are related as follows.

\begin{theorem}\label{theo-bipartite}
Let $G$ be an integrally invertible graph. Then  $G$ is bipartite if and only if $G$ is  simultaneously positively and negatively invertible.
\end{theorem}

\noindent P r o o f.
Let $G=G_B$ be an integrally invertible bipartite graph. Assume that $G_B$ is positively invertible, i.e., there exists  a signature matrix $D^+ = diag(D_1, D_2)$ such that the matrix
\begin{eqnarray*}
D^+ B^{-1} D^+ &=& 
\left( 
\begin{array}{cc}
D_1 & 0  \\
0 & D_2\\
\end{array}
\right)
\left( 
\begin{array}{cc}
0 & (K^{-1})^T  \\
K^{-1} & 0\\
\end{array}
\right)
\left( 
\begin{array}{cc}
D_1 & 0  \\
0 & D_2\\
\end{array}
\right)
\\
&=&
\left( 
\begin{array}{cc}
0 & D_1 (K^{-1})^T D_2   \\
D_2 K^{-1} D_1 & 0\\
\end{array}
\right)
\end{eqnarray*}
contains nonnegative integer entries only. Then for the $\{\pm1\}$ diagonal matrix $D^- = diag(D_1, - D_2)$ the matrix 
\[
D^- B^{-1} D^- = 
\left( 
\begin{array}{cc}
0 & -D_1 (K^{-1})^T D_2   \\
-D_2 K^{-1} D_1 & 0\\
\end{array}
\right)
\]
contains nonpositive integers only. Hence $G_B$ is negatively invertible, and vice versa.

On the other hand, suppose that $G$ is simultaneously positively and negatively invertible. We will prove that $G$ is a bipartite graph with even number of vertices. Let $n$ be the number of vertices of the graph $G$.  Let $D^\pm$ be diagonal $\{\pm 1\}$-matrices such that $D^+ A^{-1} D^+$ contains nonnegative integers and $D^- A^{-1} D^-$ contains nonpositive integers only. Since $(D^\pm A^{-1} D^\pm)_{ij} = D^\pm_{ii} (A^{-1})_{ij} D^\pm_{jj}$ we conclude that $(D^\pm A^{-1} D^\pm)_{ij}\not=0$ if and only if $(A^{-1})_{ij}\not=0$. Hence
\begin{equation}
 D^+ A^{-1} D^+ = - D^- A^{-1} D^- .
 \label{Dpm}
\end{equation}
As $det(A^{-1}) = det(D^+ A^{-1} D^+) = (-1)^n det(D^- A^{-1} D^-)= (-1)^n det(A^{-1})$ we conclude that $n$ is even, i.~e. $n=2m$. 

Recall that for the trace operator $tr(Z)=\sum_{i=1}^n Z_{ii}$ of an $n\times n$ matrix $Z$ we have $tr(X Y) = tr(Y X)$ where $X,Y$ are $n\times n$ matrices. With respect to (\ref{Dpm}) we obtain:
\begin{eqnarray*}
 tr(D^+D^-) &=& tr(A D^+ D^- A^{-1}) = tr(D^- A D^+ D^- A^{-1} D^-) 
\\
&=& - tr(D^- A D^+ D^+ A^{-1} D^+) =  -  tr(D^-D^+)=  -  tr(D^+D^-).
\end{eqnarray*}
Thus $ tr(D^+D^-)=0$. Since $D^\pm$ are diagonal $\{\pm 1\}$-matrices, there exists an $n\times n$ permutation matrix $P$ such that
\[
 P^T D^+ D^- P = 
\left( 
\begin{array}{cc}
I & 0   \\
0 & -I\\
\end{array}
\right),\ \ i.e. 
\quad 
D^+ D^- = D^- D^+ =
P \left(
\begin{array}{cc}
I & 0   \\
0 & -I\\
\end{array}
\right) P^T ,
\]
where $I$ is the $m\times m$ identity matrix. It follows from (\ref{Dpm}) that
\[
 A^{-1} = - D^+ D^-  A^{-1} D^- D^+ = 
- P \left(
\begin{array}{cc}
I & 0   \\
0 & -I\\
\end{array}
\right) 
P^T A^{-1} P
\left(
\begin{array}{cc}
I & 0   \\
0 & -I\\
\end{array}
\right) 
P^T.
\]
Since $P^T P = P P^T =I$ we have
\[
P^T A^{-1} P
= 
- \left(
\begin{array}{cc}
I & 0   \\
0 & -I\\
\end{array}
\right) 
P^T A^{-1} P
\left(
\begin{array}{cc}
I & 0   \\
0 & -I\\
\end{array}
\right) 
\]
If we write $P^T A^{-1} P$ as a block matrix we obtain 
\begin{eqnarray*}
 P^T A^{-1} P \equiv
 \left(
\begin{array}{cc}
V & H   \\
H^T & W\\
\end{array}
\right) 
&=&
- \left(
\begin{array}{cc}
I & 0   \\
0 & -I\\
\end{array}
\right) 
 \left(
\begin{array}{cc}
V & H   \\
H^T & W\\
\end{array}
\right) 
\left(
\begin{array}{cc}
I & 0   \\
0 & -I\\
\end{array}
\right) 
\\
&=&
 \left(
\begin{array}{cc}
-V & H   \\
H^T & -W\\
\end{array}
\right).
\end{eqnarray*}
Therefore $V=W=0$ and 
\[
P^T A^{-1} P = 
\left(
\begin{array}{cc}
0 & H   \\
H^T & 0\\
\end{array}
\right)
\ \Longrightarrow \ 
P^T A P = 
\left(
\begin{array}{cc}
0 & (H^T)^{-1}   \\
H^{-1} & 0\\
\end{array}
\right).
\]
This means that the adjacency matrix $A$ represents a bipartite graph $G=G_A$ after a  permutation of its vertices given by the matrix $P$. 
\hfill $\diamondsuit$

\section{Integrally invertible graphs arising by bridging}

Let $G_A$ and $G_B$ be undirected graphs on $n$ and $m$ vertices, respectively. By ${\mathcal B}_k(G_A,G_B)$ we shall denote the graph $G_C$ on $n+m$ vertices which is obtained by bridging the last $k$ vertices of the  graph $G_A$ to the first $k$ vertices of $G_B$. The adjacency matrix $C$ of the graph $G_C$ has the form:
\[
C  = \left( 
\begin{array}{cc}
A & H\\
H^T & B
\end{array}
\right),
\]
where the $n\times m$ matrix $H$ has the block structure:
\[
H= \left( 
\begin{array}{cc}
0 & 0\\
I & 0
\end{array}
\right)
= F E^T,\quad\hbox{where}\quad 
F=\left( 
\begin{array}{l}
0\\
I 
\end{array}
\right),\ \ 
E=\left( 
\begin{array}{l}
I\\
0 
\end{array}
\right),
\]
and $I$ is  the $k\times k$ identity matrix. 

Assume that $A$ and $B$ are symmetric $n\times n$ and $m\times m$ invertible matrices, respectively.  With regard to the Schur complement theorem we obtain
\begin{equation}
C^{-1}=\left( 
\begin{array}{cc}
A & H\\
H^T & B
\end{array}
\right)^{-1}
=
\left( 
\begin{array}{cc}
S^{-1} & - S^{-1} H B^{-1} \\
- B^{-1} H^T S^{-1} & B^{-1} + B^{-1} H^T S^{-1} H B^{-1}
\end{array}
\right),
\label{invC}
\end{equation}
where $S= A- H B^{-1} H^T$ is the Schur complement (see e.~g. \cite[Theorem A.6]{Maja2013}). To facilitate further notation let us introduce the following matrices:
\[
P = F^T A^{-1} F,\qquad  R = E^T B^{-1} E.
\]
In order to compute the inverse of the Schur complement $S$ we follow derivation of the Woodbury and Morrison-Sherman formulae (c.f. \cite[Corollary A.6, A.7]{Maja2013}). More precisely, equation $S x =y$ can be rewritten as follows:
\[
 y=(A- H B^{-1} H^T)x = A x - F E^T B^{-1} E F^T x = A x -F R F^T x
\]
and thus $x= A^{-1} y + A^{-1} F R F^T x$. Hence
\[
 F^T x = F^T A^{-1} y + F^T A^{-1} F R F^T x = F^T A^{-1} y + P  R F^T x .
\]
If we assume that the matrix $I - P R$ is invertible then $F^T x =  (I - P R)^{-1}F^T A^{-1} y$, and 
\begin{equation}
 S^{-1}=(A- H B^{-1} H^T)^{-1} = A^{-1} + A^{-1} F R (I - P R)^{-1}F^T A^{-1}. 
\label{invS}
\end{equation}
Note that the matrix  $I - P R$ is integrally invertible provided that either $P R=0$ or $P R = 2 I$.

Clearly, an $m\times m$  matrix with a zero principal $k\times k$ diagonal block is invertible only for $k\le m$. Consequently, there are no connected invertible graphs with $E^T B^{-1} E=0$ for $k > m/2$. 

\begin{theorem}\label{theo-1}
Let $G_A$ and $G_B$ be  integrally invertible graphs on $n$ and $m$ vertices, and let $R$ and $P$ be the upper left and lower right $k\times k$ principal submatrices of $B^{-1}$ and $A^{-1}$, respectively.

Let $G_C={\mathcal B}_k(G_A,G_B)$ be the graph obtained by bridging $G_A$ and $G_B$ over the last $k$ vertices of $G_A$ and the first $k$ vertices of $G_B$. If $P R=0$ or $P R = 2 I$, then the graph $G_C$ is integrally invertible. 
\end{theorem}

\noindent P r o o f. Since $P R=0$ or $P R = 2 I$, then the inverse $(I - P R)^{-1}$ exists and is equal to $\pm I$. Hence the inverse $S^{-1}$ of the Schur complement is an integral matrix. Therefore the block matrix $C$ given by
\[
 C = \left( 
\begin{array}{cc}
A & H\\
H^T & B
\end{array}
\right)
\]
is invertible, and hence so is the bridged graph $G_C$. Moreover, $C^{-1}$ is an integral matrix because $A^{-1}, B^{-1}, S^{-1}$ are integral. 
\hfill $\diamondsuit$

\begin{definition}\label{def-arbitrarily}
Let $G_B$ be a graph on $m$ vertices with an invertible adjacency matrix $B$. We say that $G_B$ is arbitrarily bridgeable over a subset of $k\le m/2$ vertices if the $k\times k$ upper principal submatrix $R\equiv E^T B^{-1} E$ of the inverse matrix $B^{-1}$ is a null matrix, that is $R = 0$. 
\end{definition}

In view of Definition~\ref{def-arbitrarily}, the bridged graph $G_C={\mathcal B}_k(G_A,G_B)$ is integrally invertible provided that $G_B$ is arbitrarily bridgeable over the  set of its 'first' $k$ vertices.

In the next theorem we address the question of invertibility of the bridged graph $G_C={\mathcal B}_k(G_A,G_B)$ under the assumption that $G_A$ and $G_B$ are positive (negative) invertible graphs. 

\begin{theorem}\label{theo-si}
Let $G_A$ and $G_B$ be graphs on $n$ and $m$ vertices, respectively. Assume that they are either both positively invertible or both negatively invertible graphs with signature matrices $D^A$ and $D^B$. Then the graph  $G_C={\mathcal B}_k(G_A,G_B)$ is positively (negatively) invertible if we have $P R =0$ and either the matrix $D^A H D^B$ or $-D^A H D^B$ contains nonnegative integers only.
\end{theorem}

\noindent P r o o f. Let $C$ be the adjacency matrix to the graph $G_C={\mathcal B}_k(G_A,G_B)$. If $P R =0$ then for the inverse of the Schur complement (see (\ref{invS})) we have
\[
 S^{-1} = A^{-1} + A^{-1} F R F^T A^{-1} = A^{-1} + A^{-1} H B^{-1} H^T A^{-1} 
\]
because $F R F^T = F E^T B^{-1} E F^T = H B^{-1} H^T$. Therefore
\begin{eqnarray*}
 D^A S^{-1} D^A &=& D^A A^{-1} D^A + D^A A^{-1} H B^{-1} H^T A^{-1} D^A
 \\
  &=& D^A A^{-1} D^A 
  \\
  && + (D^A A^{-1} D^A)( D^A H D^B)( D^B B^{-1} D^B)( D^B  H^T D^A)( D^A A^{-1} D^A)
\end{eqnarray*}
and so $D^A S^{-1} D^A$ is a nonnegative (nonpositive) integer matrix because the matrices $D^A A^{-1} D^A$ and $D^B B^{-1} D^B$ are simultaneously nonnegative (nonpositive) and $D^A H D^B$ or $-D^A H D^B$ contains nonnegative integers only. 

In the case when $D^A H D^B$ is nonnegative we will prove that $C^{-1}$ is diagonally similar to a nonnegative (nonpositive) integer matrix with $D^C=diag(D^A,-D^B)$\ ($D^C=diag(D^A,D^B)$). With regard to (\ref{invC}) we have
\begin{eqnarray*}
D^C C^{-1} D^C &=& 
\left( 
\begin{array}{cc}
D^A & 0\\
 0  & -D^B
\end{array}
\right)
\left( 
\begin{array}{cc}
S^{-1} & - S^{-1} H B^{-1} \\
- B^{-1} H^T S^{-1} & B^{-1} + B^{-1} H^T S^{-1} H B^{-1}
\end{array}
\right)
\\
&& \qquad\qquad \times \left( 
\begin{array}{cc}
D^A & 0\\
 0  & -D^B
\end{array}
\right)
\end{eqnarray*}
\begin{equation*}
=
\left( 
\begin{array}{cc}
D^A S^{-1} D^A  &  D^A S^{-1} D^A D^A H D^B D^B B^{-1} D^B \\
D^B  B^{-1} D^B D^B H^T D^A D^A S^{-1} D^A  & D^B B^{-1} D^B + W,
\end{array}
\right)
\end{equation*}
where 
\[
 W =  (D^B B^{-1} D^B)( D^B H^T D^A)( D^A S^{-1} D^A)( D^A H D^B)( D^B B^{-1} D^B) .
\]
In the expression for the matrix $W$ we have intentionally used the matrices $D^A D^A=D^B D^B = I$ instead of the identity matrix $I$. Since the matrices $D^A A^{-1} D^A$, $D^B B^{-1} D^B$, and  $D^A S^{-1} D^A$ contain nonnegative (nonpositive) integers only and $D^A H D^B$ is nonnegative, we conclude that $C^{-1}$ is diagonally similar to a nonnegative (nonpositive) integral matrix. 

In the case when $-D^A H D^B$ is nonnegative we can proceed similarly as before and conclude that $C^{-1}$ is diagonally similar to a nonnegative (nonpositive) integral matrix.

Hence the graph $G_C$ is positively (negatively) invertible, as claimed.
\hfill $\diamondsuit$

\bigskip

As a consequence we obtain the following:

\begin{corollary}
Let $G_A, G_B$ be two positively (negatively) invertible graphs such that $(B^{-1})_{11}=0$. Then the graph $G_C={\mathcal B}_1(G_A,G_B)$ bridged over the first vertex is again positively (negatively) invertible.
\end{corollary}
\noindent P r o o f. For $k=1$ the condition $(B^{-1})_{11}=0$ implies $R\equiv E^T B^{-1} E=0$, i.e. $G_B$ is arbitrarily bridgeable. The matrix $D^A H D^B$ contains only one nonzero element, equal to $\pm H$. Hence the assumptions of Theorem~\ref{theo-si} are fulfilled and so $G_C$ is positively (negatively) invertible.
\hfill $\diamondsuit$

\medskip

With regard to Theorem~\ref{theo-bipartite} and Theorem~\ref{theo-si} we obtain the following result:

\begin{corollary}
Let $G_A, G_B$ be two bipartite positively and negatively invertible graphs such that  $G_B$ is arbitrarily bridgeable over the first $k$ vertices. Let $D^A_+$ and $D^B_+$ ($D^A_-$ and $D^B_-$) be diagonal $\{\pm1\}$-matrices such that $D^A_+ A^{-1} D^A_+$ and  $D^B_+ B^{-1} D^B_+$ ($D^A_- A^{-1} D^A_-$ and  $D^B_- B^{-1} D^B_-$) are nonnegative (nonpositive) matrices. If $D^A_\pm  H D^B_\pm$ are either both nonnegative or both nonpositive then the bridged graph $G_C={\mathcal B}_k(G_A,G_B)$ is again bipartite positively and negatively  invertible.
\end{corollary}

\begin{figure}
\begin{center}
\includegraphics[width=7truecm]{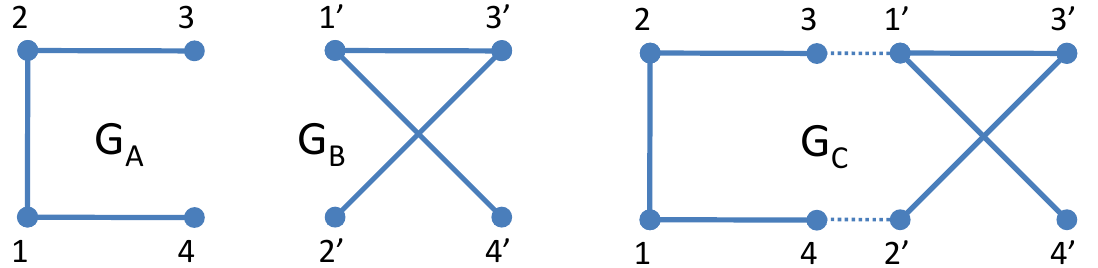}
\end{center}
\caption{
An example of bridging of two bipartite positively and negatively invertible graphs $G_A$ and $G_B$ through vertices $3\leftrightarrow 1', 4\leftrightarrow 2'$. The resulting graph $G_C={\mathcal B}_2(G_A, G_B)$. 
}
\label{fig-bridged}
\end{figure}

\begin{example}\rm 
In what follows, we present an example showing that the assumption made on nonnegativity or nonpositivity of matrices $D^A_+ H D^B_+$ and $D^A_- H D^B_-$ cannot be relaxed. To do so, we will construct a bridged graph $G_C$ from two integrally invertible bipartite graphs such that $G_C$ is only positively but not negatively invertible graph and, as a consequence of Theorem~\ref{theo-bipartite}, the graph $G_C$ is not bipartite. 

Let $G_A$ and $G_B$ be two simultaneously positively and negatively invertible bipartite graphs shown in Fig.~\ref{fig-bridged}. We will bridge them over a set of $k=2$ vertices to obtain the graph $G_C$ with inverses given by
\[
\scriptsize
A^{-1}=\left( 
\begin{array}{cccc}
0 & 0 & 0 & 1\\
0 & 0 & 1 & 0\\ 
0 & 1 & 0 & -1\\
1 & 0 & -1 & 0\\
\end{array}
\right),
\quad 
B^{-1}=\left( 
\begin{array}{cccc}
0 & 0 & 0 & 1\\
0 & 0 & 1 & -1\\ 
0 & 1 & 0 & 0\\ 
1 & -1 & 0 & 0\\
\end{array}
\right),
\]
\[\scriptsize
C^{-1}=\left( 
\begin{array}{cccccccc}
0 & 0 & 0 & 1 & 0 & 0 & -1 & 1\\
0 & 0 & 1 & 0 & 0 & 0 & 0 & -1\\
0 & 1 & 0 & -1 & 0 & 0 & 1 & -1\\
1 & 0 & -1 & 0 & 0 & 0 & 0 & 1\\
0 & 0 & 0 & 0 & 0 & 0 & 0 & 1\\
0 & 0 & 0 & 0 & 0 & 0 & 1 & -1\\
-1 & 0 & 1 & 0 & 0 & 1 & 0 & -1\\
1 & -1 & -1 & 1 & 1 & -1 & -1 & 2\\
\end{array}
\right).
\]
The graphs $G_A, G_B$ are isomorphic with eigenvalues: $\{\pm 1.6180, \pm 0.6180\}$. 
The upper left $2\times 2$ principal submatrix $R$ of $B^{-1}$ is zero, so that the graph $G_B$ can be arbitrarily bridged to an integrally invertible graph $G_A$. 

It is easy to verify that the inverse matrices $A^{-1}$ and $B^{-1}$ can be signed to nonnegative matrices by signature matrices $D^A_+=D^B_+ = diag(-1,1,1,-1)$. At the same time they can be signed to a nonpositive matrix by $D^A_-=diag(-1,-1,1,1)$ and $D^B_-=diag(-1,1,-1,1)$. Furthermore, $D^A_+ H D^B_+ = -H$ is a nonpositive matrix. By Theorem~\ref{theo-si}, the graph $G_C$ is positively invertible. On the other hand, 
\[\scriptsize 
D^A_- H D^B_- = \left( 
\begin{array}{cccc}
0 & 0 & 0 & 0\\
0 & 0 & 0 & 0\\ 
-1 & 0 & 0 & 0\\ 
0 & 1 & 0 & 0\\
\end{array}
\right)
\]
is neither nonnegative nor nonpositive. Indeed, the graph $G_C$ is not bipartite and it is just positively (and not negatively) invertible, with spectrum
\[
\sigma(G_C)= \{-1.9738, -1.8019, -0.7764, -0.445, 0.2163, 1.247, 1.4427, 2.0912\}.
\] 
\end{example}

\begin{remark}
If the graph $G_B$ is arbitrarily bridgeable over the first $k$ vertices then $R=0$, and, consequently the assumption $P R =0$ appearing in Theorem~\ref{theo-si} is satisfied. On the other hand, if we consider the graph $G_C$ with the vertex set $\{1,2,3,4,1',2',3',4'\}$ shown in Figure~\ref{fig-bridged} then the inverse matrix $C^{-1}$ contains the principal submatrix 
\[\scriptsize
 P = \left( 
\begin{array}{cccc}
0 & 0 & 0 & 0\\
0 & 0 & 0 & 0\\ 
0 & 0 & 0 & 1\\ 
0 & 0 & 1 & 0\\
\end{array}
\right)
\]
corresponding to vertices $4,1',2',3'$. Consider the same graph $\tilde{G}_{\tilde{C}}$ on the vertex set $\{\tilde{1},\tilde{2},\tilde{3},\tilde{4},\tilde{1}',\tilde{2}',\tilde{3}',\tilde{4}'\}$. Then, after permuting vertices, the inverse matrix $\tilde{C}^{-1}$ has the upper left principal $4\times 4$ submatrix 
\[\scriptsize
  R = \left( 
\begin{array}{cccc}
0 & 1 & 0 & 0\\
1 & 0 & 0 & 0\\ 
0 & 0 & 0 & 0\\ 
0 & 0 & 0 & 0\\
\end{array}
\right).
\]
Hence $P R =0$ but neither $P$ nor $R$ is an all-zero $4\times4$ matrix. By Theorem~\ref{theo-si} the bridged graph ${\mathcal B}_4(G_C, \tilde{G}_{\tilde{C}})$ over the set of vertices $4\leftrightarrow\tilde{2}', 1'\leftrightarrow \tilde{3}', 2'\leftrightarrow \tilde{4}, 3'\leftrightarrow \tilde{1}'$ is integrally invertible.

\end{remark}

\section{Spectral bounds for graphs arising by bridging}

In this section we derive a lower bound for the least positive eigenvalue of  bridged graphs  ${\mathcal B}_k(G_A,G_B)$ in terms of the least positive eigenvalues of graphs $G_A$ and $G_B$. Throughout this  section we assume that the adjacency matrices $A,B$ are invertible but we do not require their integral invertibility. 

Before stating and proving our spectral estimate we need the following auxiliary Lemma.

\begin{lemma}\label{lemmaMax}
Assume that $D$ is an $n\times m$ matrix and $\alpha,\beta>0$ are positive constants. Then, for the optimal value $\lambda^*$ of the following constrained optimization problem:
\begin{equation}\label{minD}
\begin{array}{rl}
\lambda^*= \max & \alpha \Vert x- D y\Vert^2 +\beta \Vert y\Vert^2 \\
{\rm s. t.} & \Vert x\Vert^2 + \Vert y\Vert^2  = 1,\ \  x\in\R^n, y\in\R^m,
\end{array} 
\end{equation}
we have an explicit expression of the form:
\[
 \lambda^* = \max\left\{ \lambda,\ \frac{(\lambda-\alpha)(\lambda-\beta)}{\alpha\lambda} \in\sigma(D^T D)\right\} 
\]
\[
 =  \frac{\alpha\mu^* + \alpha+\beta + \sqrt{(\alpha\mu^* + \alpha+\beta)^2-4\alpha\beta}}{2},
\]
where $\mu^*=\max\{\sigma(D^T D)\}$ is the maximal eigenvalue of the matrix $D^T D$. 
\end{lemma}

\noindent P r o o f.
The proof is straightforward and is based on standard application of the Lagrange multiplier method (see e.g. \cite{Maja2013}). We give details for the reader's convenience.

Let us introduce the Lagrange function:
\begin{eqnarray*}
 L(x,y,\lambda) &=& \alpha \Vert x- D y\Vert^2 +\beta \Vert y\Vert^2 - \lambda( \Vert x\Vert^2 + \Vert y\Vert^2) 
 \\
 &=& \alpha x^T x -2\alpha x^T D y +\alpha y^T D^T D y +\beta y^T y -\lambda x^T x -\lambda y^T y.
\end{eqnarray*}
Now, it follows from the first order necessary conditions for constrained maximum $(x,y)$ (see e.~g. \cite{Maja2013}) that there exists a Lagrange multiplier $\lambda\in\R$ such that 
\begin{eqnarray}
\label{A1}
0 &=& L^\prime_x \equiv 2 \alpha x^T  -2\alpha y^T D^T -2 \lambda x^T,
\\
 0 &=& L^\prime_y \equiv  -2\alpha x^T D  +2 \alpha y^T D^T D + 2 \beta y^T -2 \lambda y^T.
 \end{eqnarray}
In the case $\lambda\not=\alpha$ we obtain
\[
x = \frac{\alpha}{\alpha-\lambda} D y,
\qquad  
\left[ \alpha D^T D -(\lambda-\beta)\right] y = \alpha D^T x =  \frac{\alpha^2}{\alpha-\lambda}  D^T Dy. 
\]
Therefore, 
\[
  D^T D y = \frac{(\lambda-\alpha)(\lambda-\beta)}{\alpha\lambda} y \ \Rightarrow \ 
  \frac{(\lambda-\alpha)(\lambda-\beta)}{\alpha\lambda} \in\sigma(D^T D).
\]
Now, from the constraint $x^T x + y^T y =1$ we deduce that 
\[
 1 = x^T x + y^T y =  \frac{\alpha^2}{(\alpha-\lambda)^2} y^T D^T D y  + y^T y = 
 \left( \frac{\alpha^2}{(\alpha-\lambda)^2}  \frac{(\lambda-\alpha)(\lambda-\beta)}{\alpha\lambda} +1 \right) \Vert y\Vert^2.
 \]
Hence
\[
 \Vert y\Vert^2 =  \frac{(\lambda-\alpha)\lambda}{\lambda^2 -\alpha\beta}, \qquad
 \Vert x\Vert^2 = 1- \Vert y\Vert^2 =  \frac{(\lambda-\beta)\alpha}{\lambda^2 -\alpha\beta}.
\]
Finally, for the value function $f(x,y) = \alpha \Vert x- D y\Vert^2 +\beta \Vert y\Vert^2 $ of the constrained optimization problem (\ref{minD}) we obtain
\begin{eqnarray*}
 f(x,y) &=&  \alpha x^T x -2\alpha x^T D y +\alpha y^T D^T D y +\beta y^T y
 \\
 &=&  \alpha x^T x -2\frac{\alpha^2}{\alpha-\lambda} y^T D^T D y +\alpha y^T D^T D y +\beta y^T y
 \\
 &=& \frac{\alpha^2(\lambda-\beta)}{\lambda^2-\alpha\beta}
 + \left(\alpha - \frac{2\alpha^2}{\alpha-\lambda}\right)\frac{(\lambda-\alpha)(\lambda-\beta)}{\alpha\lambda} \frac{(\lambda-\alpha)\lambda}{\lambda^2 -\alpha\beta} 
 +
 \beta \frac{(\lambda-\alpha)\lambda}{\lambda^2 -\alpha\beta} 
 \\
 &=& \lambda.
\end{eqnarray*}
In the case $\lambda=\alpha$, one sees from (\ref{A1}) that $D y=0$ and so $f(x,y)= \alpha \Vert x\Vert^2 +\beta \Vert y\Vert^2 \le \max\{\alpha,\beta\} \le \lambda^*$. In summary, 
\[\lambda^* = \max\left\{ \lambda,\ \frac{(\lambda-\alpha)(\lambda-\beta)}{\alpha\lambda} \in\sigma(D^T D)\right\}, 
\]
as claimed. \hfill$\diamondsuit$

We are in a position to present our spectral bound.

\begin{theorem}\label{theo-2}
Let $G_A$ and $G_B$ be  graphs on $n$ and $m$ vertices with invertible adjacency matrices. Assume $G_B$ is arbitrarily bridgeable over the first $k$ vertices.  Then the least positive eigenvalue $\lambda_1^+(G_C)$ of its adjacency matrix $C$ of the bridged graph $G_C={\mathcal B}_k(G_A,G_B)$ satisfies 
\[
 \lambda_1^+(G_C) \ge \lambda_{lb}(G_A,G_B,k) :=\frac{2}{\alpha\mu^* + \alpha+\beta + \sqrt{(\alpha\mu^* + \alpha+\beta)^2-4\alpha\beta}},
\]
where $\mu^*=\max\{\sigma(B^{-1} H^T H B^{-1})\}$ is the maximal eigenvalue of the positive semidefinite $m\times m$ matrix $B^{-1} H^T H B^{-1}, \alpha = 1/\lambda_1^+(G_A)$ and $\beta = 1/\lambda_1^+(G_B)$.
\end{theorem}

\noindent P r o o f.

The idea of the proof is based on estimation of the numerical range of the matrix $C^{-1}$. Since 
$\lambda_1^+(G_C)=\lambda_1^+(C)=1/\lambda_{max}(C^{-1})$ where $\lambda_{max}(C^{-1})$ is the maximal eigenvalue of the inverse matrix $C^{-1}$,  the lower bound for $\lambda_1^+(C)$ can be derived from the upper bound for $\lambda_{max}(C^{-1})$. As stated in Definition~\ref{def-arbitrarily}, the assumption that $G_B$ is an arbitrarily bridgeable graph implies $S^{-1}=A^{-1}$. Thus formula (\ref{invC}) for the inverse matrix $C^{-1}$ becomes:

\begin{eqnarray*}
C^{-1} &=&
\left( 
\begin{array}{cc}
A^{-1} & - A^{-1} H B^{-1} \\
- B^{-1} H^T A^{-1} & B^{-1} + B^{-1} H^T A^{-1} H B^{-1}
\end{array}
\right)
\\
&=& 
\left( 
\begin{array}{cc}
I & 0 \\
- B^{-1} H^T & I
\end{array}
\right)
\left( 
\begin{array}{cc}
A^{-1} & 0 \\
0  & B^{-1}
\end{array}
\right)
\left( 
\begin{array}{cc}
I & - H B^{-1} \\
0  & I
\end{array}
\right).
\end{eqnarray*}
Let $z=(x,y)^T\in \R^{n+m}$ where $x\in\R^n, y\in\R^m$. 
For the Euclidean inner product $\langle C^{-1} z, z \rangle $ in $\R^{n+m}$ we obtain
\begin{eqnarray*}
\langle C^{-1} z, z\rangle  &=& \langle A^{-1} (x-H B^{-1} y), x-H B^{-1} y \rangle 
+  \langle B^{-1} y, y \rangle 
\\
&\le& \lambda_{max}(A^{-1}) \Vert x-H B^{-1} y\Vert^2 + \lambda_{max}(B^{-1}) \Vert y\Vert^2 .
\end{eqnarray*}
Letting $\alpha=\lambda_{max}(A^{-1}), \beta=\lambda_{max}(B^{-1})$ and $D=H B^{-1}$, by Lemma~\ref{lemmaMax} we obtain
\[
 \langle C^{-1} z, z\rangle \le \frac{1}{\lambda_{lb}(G_A,G_B,k)} \Vert z\Vert^2, 
\]
for any $z\in\R^{n+m}$. Since
\[
\lambda_{max}(C^{-1}) = \max_{z\not=0} \frac{ \langle C^{-1} z, z\rangle}{\Vert z\Vert^2} 
\le  \frac{1}{\lambda_{lb}(G_A,G_B,k)}
\]
our Theorem follows because $\alpha = \lambda_{max}(A^{-1})= 1/\lambda_1^+(A)= 1/\lambda_1^+(G_A)$ and $\beta = \lambda_{max}(B^{-1})= 1/\lambda_1^+(B)= 1/\lambda_1^+(G_B)$.
\hfill$\diamondsuit$

\bigskip

To illustrate this on an example, for the graph $G_C = {\mathcal B}_k(G_A,G_B)$ shown in Fig.~\ref{fig-bridged} we have $\lambda_1^+(G_C) = 0.2163$ and the lower bound derived above  gives $\lambda_{lb}(G_A,G_B,k) = 0.1408$.

\section{A ``fulvene'' family of integrally invertible graphs}

The aim of this section is to present construction of a family of integrally invertible graphs grown from the ``fulvene'' graph of Fig.~\ref{fig-fulvene} (left), which is the same as the $H_{10}$ in Fig.~\ref{fig-hexafamily}.  With regard to Table~\ref{tab-2} (see Section~6), the graph $F_0\equiv H_{10}$ can be arbitrarily bridged over the pair of vertices labeled by $1,2$ (see the left part of Fig.~\ref{fig-fulvene}) to any integrally invertible graph. 

\begin{figure}
\begin{center}
\includegraphics[angle=90,width=1.2truecm]{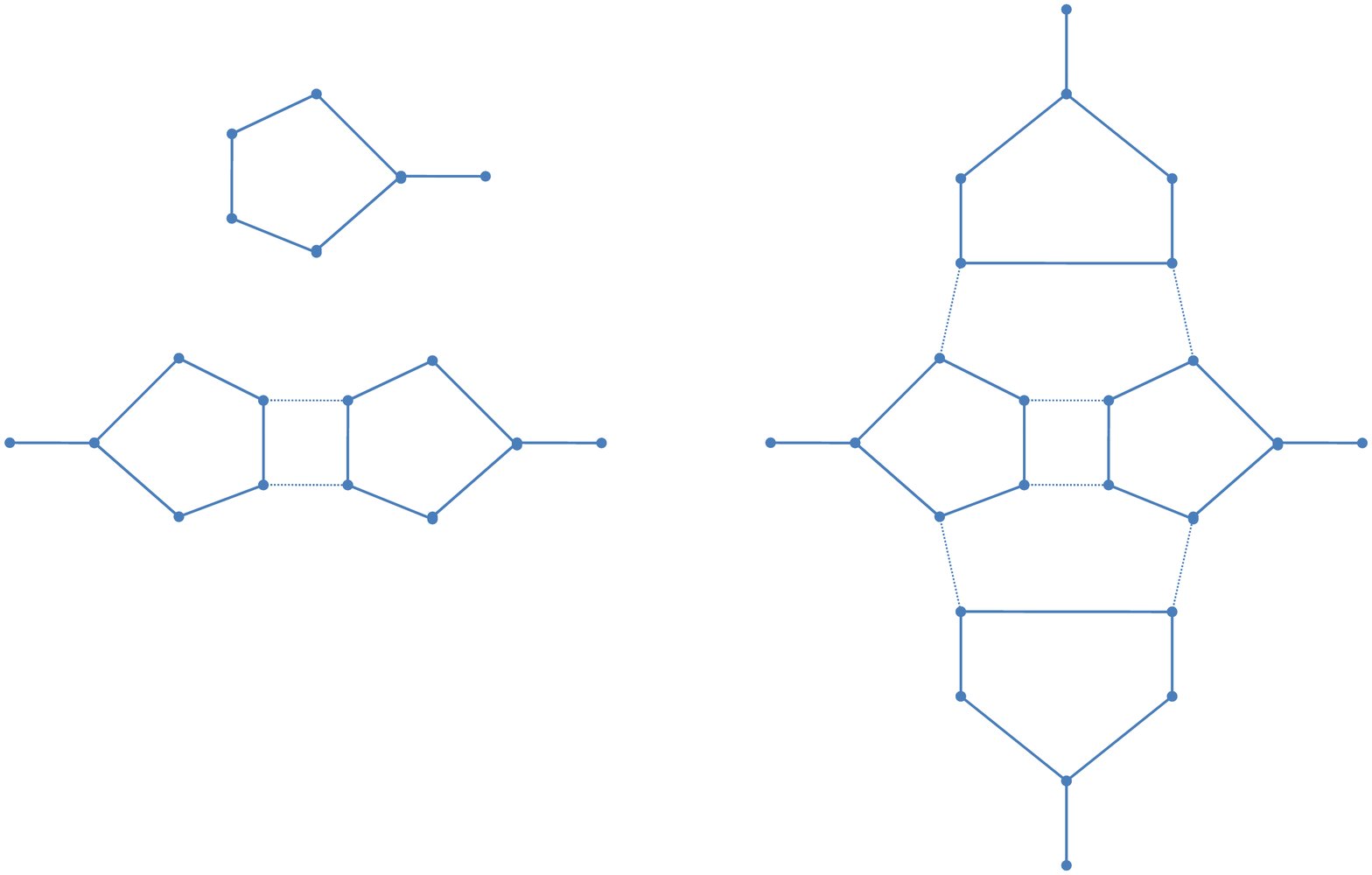}
\qquad 
\includegraphics[angle=90,width=1truecm]{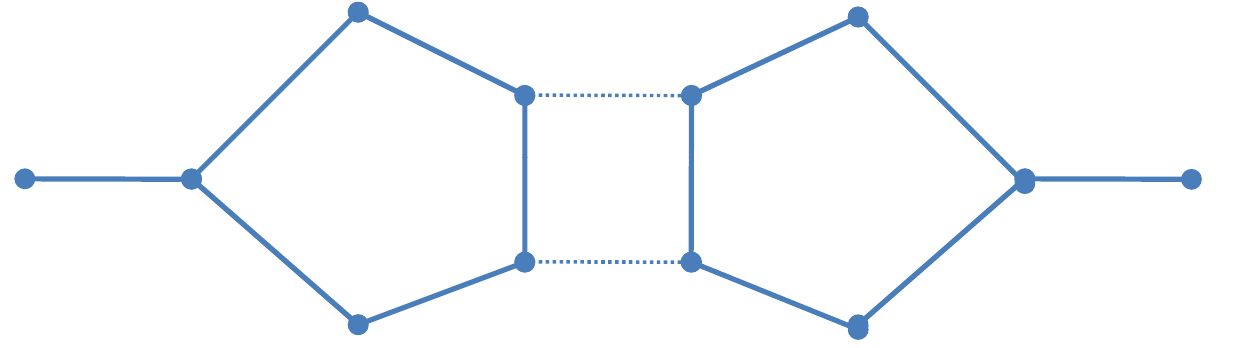}
\qquad
\includegraphics[angle=90,width=5truecm]{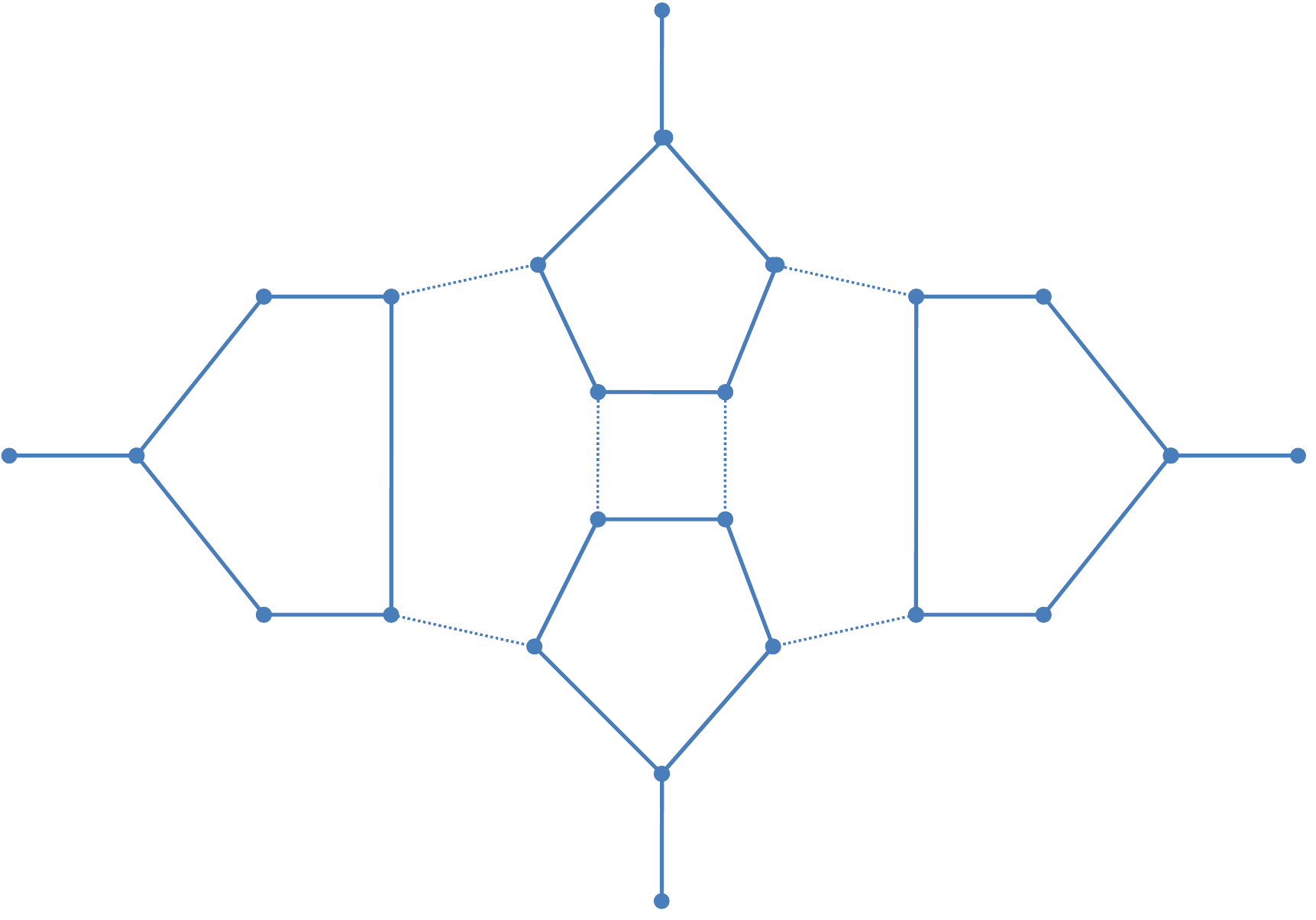}
\\
\ \hglue-2truecm $F_0$ \hskip 2truecm $F_1$ \hskip 3truecm $F_2$

\includegraphics[width=6.5truecm]{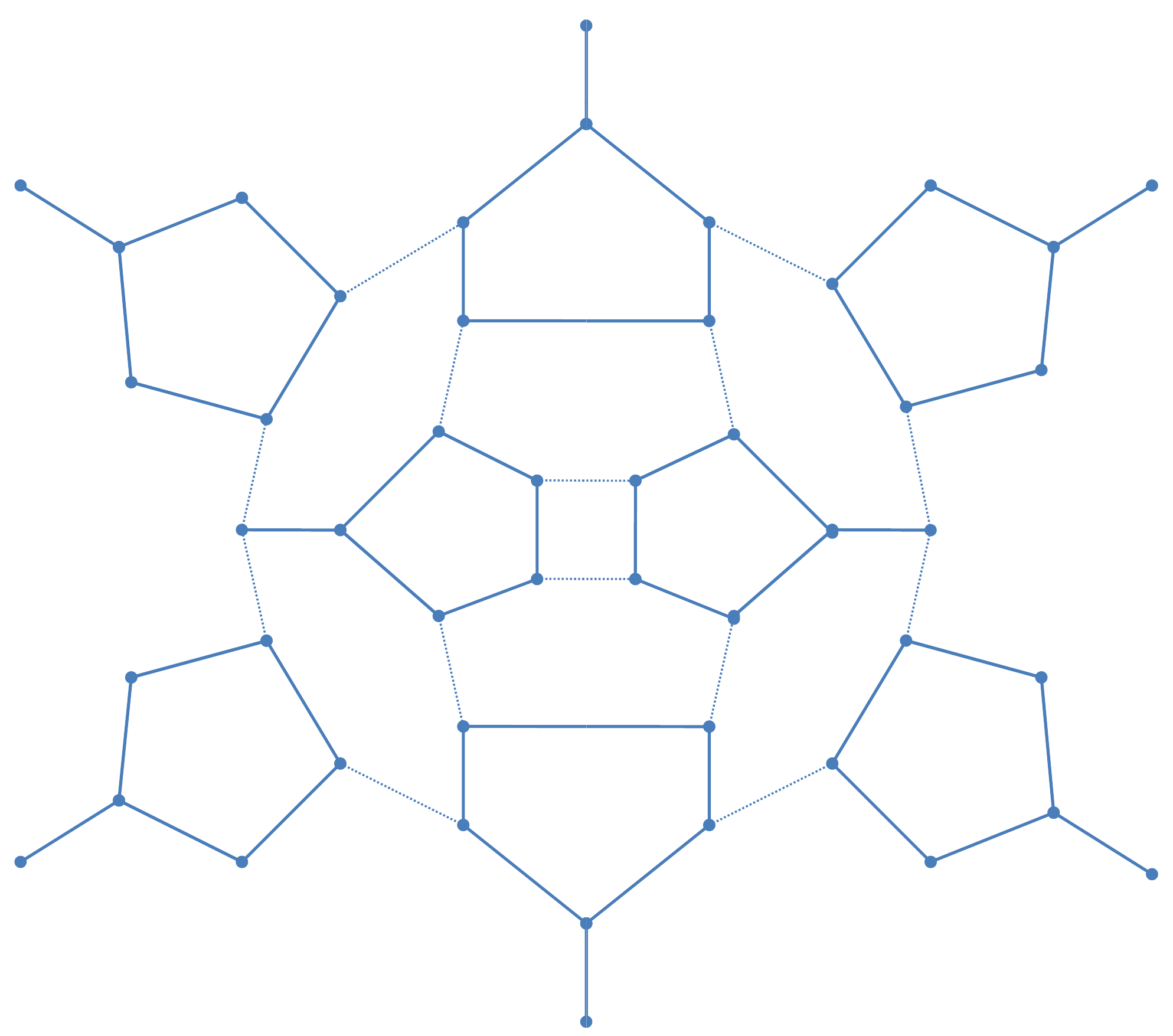}
\\
$F_3$
\end{center}
\caption{
The graphs $F_0,F_1,F_2, F_3$ of the ``fulvene'' family of integrally invertible graphs.
}
\label{fig-fulvenefamily1}
\end{figure}

Our construction begins with the fulvene graph $F_0$. The next iteration $F_1$ is obtained by bridging $F_0$ to another copy of $F_0$ over the vertex set  $\{1,2\}$ in both copies (see Fig.~\ref{fig-fulvenefamily1}).

We now describe a recursive construction of graphs $F_n$ from $F_{n-1}$. For $n\ge 2$, the graph $F_n$ will be obtained from $F_{n-1}$ by bridging a certain number $f_n$ (to be described below) copies of the graph $F_0$ over the vertex set $\{1,2\}$ to vertices of $F_{n-1}$ of degree $1$ or $2$. By definition, we set $f_1=f_2:=2$. Let $|V^{(i)}(F_n)|,\ i=1,2,3,$  denote the number of vertices of $F_n$ with degree $i$.
 
The order of bridging is as follows:
\begin{itemize}
\item two copies of $F_0$ are bridged to every vertex of degree $1$ which belonged to  $F_{n-2}$ and remained in $F_{n-1}$ with degree 1. The other vertex of $F_0$ is bridged to the shortest path distance vertex of degree $2$ belonging to $F_{n-1}$. The number $f^{(1)}_n$ of graphs $F_0$ added to $F_{n-1}$ is given by: $f^{(1)}_n = 2 f_{n-2}$. This way one uses $|V^{(2)}(F_{n-1})| -  2 f_{n-2}$ of vertices of degree $2$ from $F_{n-1}$.
 
\item The remaining $f^{(2)}_n = f_n -f^{(1)}_n$ copies of $F_0$ are bridged to $F_{n-1}$ through $|V^{(2)}(F_{n-1})| -  2 f_{n-2}$ vertices of degree $2$ in such a way that the graph is bridged to the pair vertices of degree $2$ with the shortest distance. By construction we have $|V^{(2)}(F_n)|= 2 f_n$. Hence, the number $|V^{(2)}(F_{n-1})| -  2 f_{n-2} = 2( f_{n-1} - f_{n-2})$ is even and so $f^{(2)}_n = f_{n-1} - f_{n-2}$. Moreover, the number $|V^{(1)}(F_n)|$ of vertices of degree $1$ is given by: $|V^{(1)}(F_n)| = f_n + f_{n-1}$ as the vertices of degree $1$ from $F_{n-1}\setminus F_{n-2}$ have not been bridged.
 
\end{itemize}

Since $f^{(1)}_n = 2 f_{n-2}$ and $f^{(2)}_n = f_{n-1} - f_{n-2}$ the total number $f_n=f^{(1)}_n + f^{(2)}_n$ of newly added  graphs $F_0$ satisfies the Fibonacci recurrence
\[
f_n=f_{n-1} + f_{n-2}, \ f_1=f_2=2.
\]
It can be explicitly expressed as: 
\[
 f_n = \frac{2}{\sqrt{5}} \left( q^n - q^{-n}\right),
\]
where $q=(1+\sqrt{5})/2$ is the golden ratio.

Properties of the fulvene family of graphs $F_n$, $n\ge 0$,  (for $F_0,F_1,F_2,F_3$ see Fig.~\ref{fig-fulvenefamily1}) can be summarized as follows.

\begin{theorem}\label{theo-fulvene}
Let $F_n, n\ge 0$, be a graph from the fulvene family of graphs. Then
\begin{enumerate}
\item $F_n$ is integrally invertible;
\item $F_n$ is a planar graph of maximum degree 3, with
\begin{eqnarray*}
|V^{(1)}(F_n)| &=& f_n + f_{n-1}, \\
|V^{(2)}(F_n)| &=& 2 f_n, \\
|V^{(3)}(F_n)| &=& |V(F_n)| -  |V^{(1)}(F_n)| - |V^{(2)}(F_n)| = 6 \sum_{k=1}^n f_k - 3 f_n - f_{n-1},
\end{eqnarray*}
where $|V(F_n)| = 6 \sum_{k=1}^n f_k$ is the number of vertices of $F_n$;

\item $F_n$ is asymptotically cubic in the sense that
\[
 \lim_{n\to\infty} \frac{|V^{(3)}(F_n)|}{|V(F_n)|} =1;
\]
\item the least positive eigenvalue $\lambda_1^+(F_n)$ satisfies the estimate:
\[
\lambda_1^+(F_n)\ge \frac{1}{q} \frac{5}{6^{n+1}-1}.
\]
\end{enumerate}

\end{theorem}

\noindent P r o o f. The number of vertices and integral invertibility of $F_n$ have been derived during construction of $F_n$. 

To prove the lower bound for the least positive eigenvalue $\lambda_1^+(F_n)$ of the integrally invertible  graph $F_n$ constructed in Section 5. Recall that the next generation $F_n$ is constructed from $F_{n-1}$ by bridging $f_n$ basic fulvene graphs $F_0$ to vertices of degree 1 and 2 of $F_{n-1}$, which can be described as 
\[
 F_n = {\mathcal B}_{2f_n} (F_{n-1}, G_{B_n}),
\]
where the graph $G_{B_n}$ has an $M\times M$ adjacency matrix $B_n$ of the block diagonal form:
\[
 B_n =diag( \underbrace{B,\cdots,B}_{f_n\ times}).
\]
Here $M=6f_n$ and $B=A_{F_0}$ is the adjacency matrix to the graph $F_0$. Therefore
\[
 A_{F_n} =
\left(\begin{array}{cc}
A_{F_{n-1}} & H_n\\
H_n^T & B_n
 \end{array}
 \right)
\]
where $H_n = (H_n^1, \cdots, H_n^{f_n})$ is an $N\times M$ block matrix with $N=|V(F_{n-1})|$. Each $H_n^r$ is an $N\times 6$ $\{0,1\}$-matrix of the form $H_n^r=(u^r, v^r, 0, 0, 0, 0)$ where $u^r_i=1$ ($v^r_i=1$) if and only if the vertex 1 (2) of the $r$-th fulvene graph $F_0$ is bridged to the $i$-th vertex of $F_{n-1}$. 

In order to apply the spectral estimate from Theorem~\ref{theo-2} we will derive an upper bound on the optimum value of $\mu^* = \max\sigma(B_n^{-1} H_n^T H_n B_n^{-1})$. Clearly, the matrix $B_n^{-1} H_n^T H_n B_n^{-1}$ satisfies
\[
(B_n^{-1} H_n^T H_n B_n^{-1})_{rs} = B^{-1} (H_n^r)^T H_n^s B^{-1}.
\]
Now,
\[
(H_n^r)^T H_n^s =
 \left\{
 \begin{array}{ll}
diag(1,1,0,0,0,0), & \hbox{if } r=s,\\
\\
diag(1,0,0,0,0,0), & \hbox{if $r\not=s$ and the $r$-th and $s$-th graph $F_0$} \\
                  & \hbox{are bridged to the same vertex of $F_{n-1}$,} \\
                  \\
0, & \hbox{otherwise.} 
 \end{array}
 \right.
\]
Since
\[\scriptsize 
B^{-1}
=\left(
\begin{array}{rrrrrr}
0 & 0 & 0 & 0 & 1 & -1\\
0 & 0 & 1 & 0 & 0 & -1\\
0 & 1 & 0 & 0 & -1 & 1\\
0 & 0 & 0 & 0 & 0 & 1\\
1 & 0 & -1 & 0 & 0 & 1\\
-1 & -1 & 1 & 1 & 1 & -2
 \end{array}
 \right)
\]
it can be verified by an easy calculation that
\[
\max\sigma(B^{-1} (H_n^r)^T H_n^s B^{-1}) =
 \left\{
 \begin{array}{ll}
3, & \hbox{if } r=s,\\
\\
2, & \hbox{if $r\not=s$ and the $r$-th and $s$-th graph $F_0$} \\
                  & \hbox{are bridged to the same vertex of $F_{n-1}$,} \\
                  \\
0, & \hbox{otherwise.} 
 \end{array}
 \right.
\]
Thus, for any vector $z=(z^1, \cdots, z^{f_n}) \in \R^M, z^i\in\R^6$, we have 
\begin{eqnarray*}
z^T  B_n^{-1} H_n^T H_n B_n^{-1} z 
&=& \sum_{r,s=1}^{f_n} (z^r)^T B^{-1} (H_n^r)^T H_n^s B^{-1} z^s
\\
&\le& 3\Vert z\Vert^2 + \sum_{r\not=s} (z^r)^T B^{-1} (H_n^r)^T H_n^s B^{-1} z^s
\\
&\le& 3\Vert z\Vert^2 + 2 \sum_{r\not=s} \frac12 (\Vert z^r\Vert^2 + \Vert z^s\Vert^2)\le 5\Vert z\Vert^2,
\end{eqnarray*}
because for the symmetric matrix $W=B^{-1} (H_n^r)^T H_n^s B^{-1}$ it holds:
$|u^T W v| \le \max |\sigma(W)| \frac12 (\Vert u\Vert^2 + \Vert v\Vert^2)$. Hence,
\[
 \mu^* = \max\sigma(B_n^{-1} H_n^T H_n B_n^{-1}) \le 5. 
\]
Finally, we establish the lower bound for the least positive eigenvalue $\lambda_1^+(F_n)$. With regard to Theorem~\ref{theo-2} we have
\[
\lambda_1^+(F_n)\ge 
\frac{2}{\alpha(\mu^* + 1)+\beta + \sqrt{(\alpha(\mu^* +1)+\beta)^2-4\alpha\beta}},
\]
where $\alpha=1/\lambda_1^+(F_{n-1}), \beta=1/\lambda_1^+(G_{B_n}) = 1/\lambda_1^+(F_0)=q$. If we denote $y_n=1/\lambda_1^+(F_n)$ we obtain
\begin{eqnarray*}
 y_n &\le& \frac12 \left( (\mu^*+1) y_{n-1} + q + \sqrt{ ((\mu^*+1) y_{n-1} +q)^2 -4 q y_{n-1} }\right)\\
 &\le& (\mu^*+1) y_{n-1} + q \le 6 y_{n-1} + q.
\end{eqnarray*}
Solving the above difference inequality yields $y_n\le \frac{q}{5}( 6^{n+1}-1)$ and so
\[
\lambda_1^+(F_n)\ge \frac{1}{q} \frac{5}{6^{n+1}-1},
\]
as claimed. \hfill $\diamondsuit$

\medskip

\begin{remark}
The asymptotic behavior of $\lambda_1^+(F_n)\to0$ as $n\to\infty$ is not surprising. For example, if we consider a cycle $C_N$ on $N$ vertices then, we have $\lambda_1^+(C_N)= 2\cos\frac\pi2 \frac{N-1}{N}$ and so  $\lambda_1^+(C_N)= O(N^{-1})=O(|V(C_N)|^{-a})$ with the polynomial decay rate $a=1$. In the case of the graph $F_n$ the number of its vertices growths exponentially $N=O(q^{n+1})$, and so the  lower bound $\lambda_1^+(F_n) \ge  O(6^{-n-1}) = O(|V(F_n)|^{-a})$ with the polynomial decay rate $a= \ln 6/\ln q \dot{=} 3.7234$ as $|V(F_n)|\to\infty$ can be expected.

\end{remark}

\section{Arbitrarily bridgeable connected graphs with a unique 1-factor}

In this section we present a census of invertible graphs on $m\le 6$ vertices with a unique 1-factor, such that they can be arbitrarily bridged to an invertible graph through a set of $k\le m/2$ vertices. Recall that a graph $G$ has a unique 1-factor if $G$ contains a unique 1-regular spanning subgraph (i.e., a perfect matching). Note that any graph having a 1-factor should have even number of vertices. 

For $m=2$ the graph $K_2$ is the unique connected graph with  a unique 1-factor. It is a positively invertible bipartite graph with the spectrum $\sigma(K_2)=\{-1,1\}$. 

\begin{figure}
\begin{center}
\includegraphics[width=0.5\textwidth]{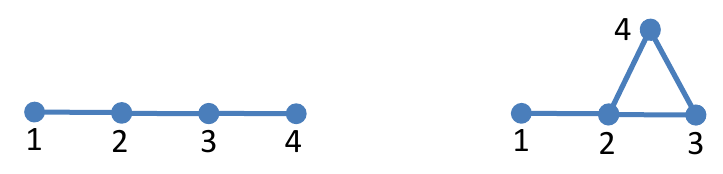}
\\
$Q_1$ \hskip  0.3\textwidth $Q_2$
\end{center}
\caption{
The family of graphs on $4$ vertices with a unique 1-factor.
}
\label{fig-quarticfamily}
\end{figure}

For $m=4$ there are two connected graphs $Q_1, Q_2$ with a unique 1-factor shown in  Fig.~\ref{fig-quarticfamily}. Both graphs are positively invertible with the spectra 
\[
\sigma(Q_1) = \{\pm 1.6180, \pm 0.6180 \},\qquad \sigma(Q_2)=\{-1.4812, -1, 0.3111, 2.1701\}.
\] 
The graph $Q_1$ can be arbitrarily bridged over the singleton sets  $\{1\}$, $\{2\}$, $\{3\}$, $\{4\}$  and over  pairs of  vertices: $\{2,3\}, \{1,3\}, \{2,4\}$. The graph $Q_1$ can be arbitrarily bridged over the singletons $\{2\}, \{3\}, \{4\}$ and over the pair $\{2,3\}$.

\begin{figure}
\begin{center}
\includegraphics[width=0.3\textwidth]{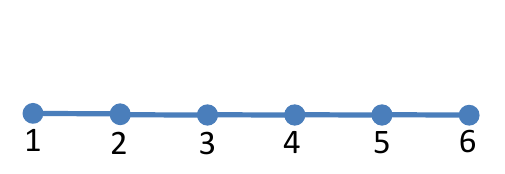}
\quad 
\includegraphics[width=0.25\textwidth]{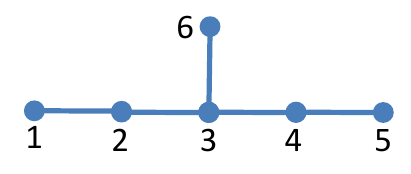}
\quad 
\includegraphics[width=0.25\textwidth]{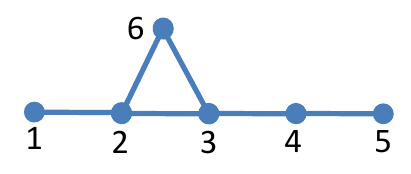}
\\
$H_1$ \hskip  0.3\textwidth $H_2$  \hskip  0.3\textwidth $H_3$
\smallskip

\includegraphics[width=0.24\textwidth]{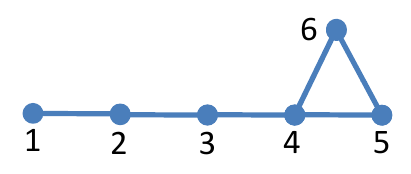}
\hskip  0.08\textwidth
\includegraphics[width=0.18\textwidth]{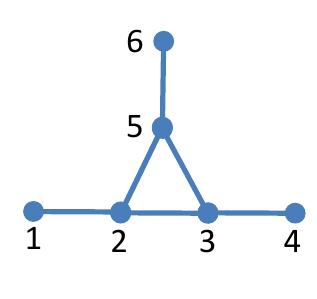}
\hskip  0.15\textwidth
\includegraphics[width=0.12\textwidth]{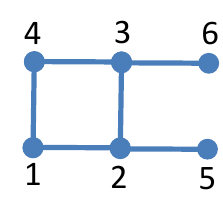}
\\
$H_4$ \hskip  0.3\textwidth $H_5$  \hskip  0.3\textwidth $H_6$
\medskip

\includegraphics[width=0.12\textwidth]{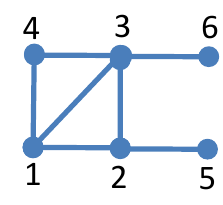}
\hskip  0.19\textwidth
\includegraphics[width=0.18\textwidth]{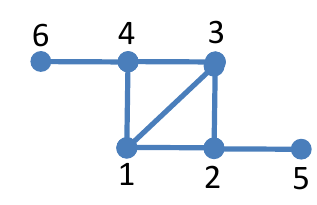}
\hskip  0.19\textwidth
\includegraphics[width=0.12\textwidth]{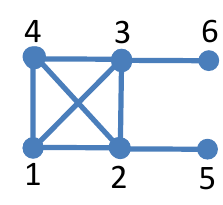}
\\
$H_7$ \hskip  0.3\textwidth $H_8$  \hskip  0.3\textwidth $H_9$
\medskip

\includegraphics[width=0.16\textwidth]{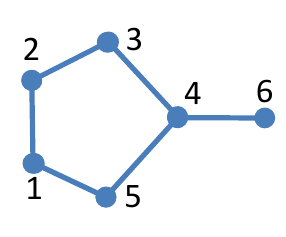}
\hskip  0.19\textwidth
\includegraphics[width=0.16\textwidth]{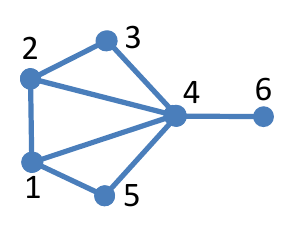}
\hskip  0.19\textwidth
\includegraphics[width=0.16\textwidth]{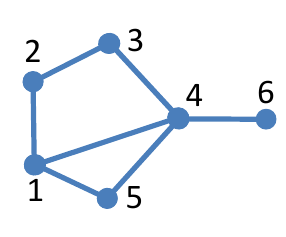}
\\
$H_{10}$ \hskip  0.3\textwidth $H_{11}$  \hskip  0.3\textwidth $H_{12}$
\medskip

\includegraphics[width=0.18\textwidth]{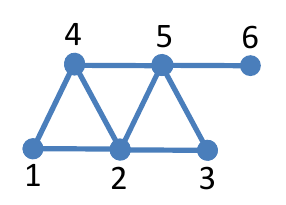}
\hskip  0.19\textwidth
\includegraphics[width=0.12\textwidth]{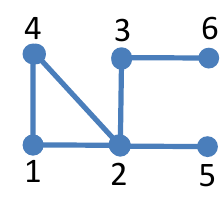}
\hskip  0.19\textwidth
\includegraphics[width=0.18\textwidth]{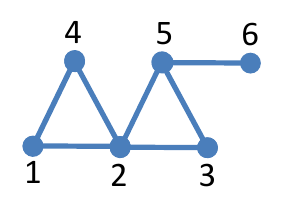}
\\
$H_{13}$ \hskip  0.3\textwidth $H_{14}$  \hskip  0.3\textwidth $H_{15}$
\medskip

\includegraphics[width=0.18\textwidth]{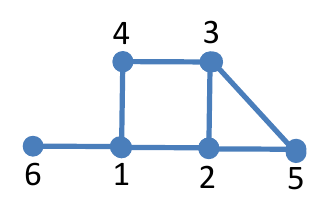}
\hskip  0.19\textwidth
\includegraphics[width=0.12\textwidth]{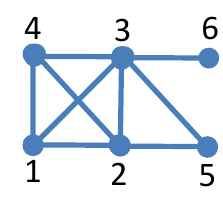}
\hskip  0.19\textwidth
\includegraphics[width=0.12\textwidth]{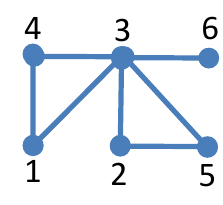}
\\
$H_{16}$ \hskip  0.3\textwidth $H_{17}$  \hskip  0.3\textwidth $H_{18}$
\medskip

\includegraphics[width=0.18\textwidth]{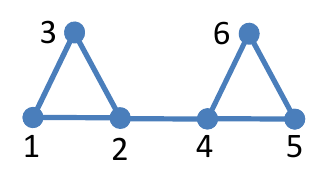}
\hskip  0.19\textwidth
\includegraphics[width=0.12\textwidth]{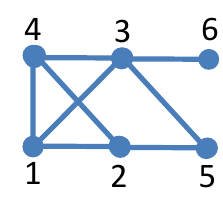}
\\
$H_{19}$\hskip  0.3\textwidth $H_{20}$ 
\medskip
\end{center}
\caption{
The family of graphs on $6$ vertices with a unique 1-factor.
}
\label{fig-hexafamily}
\end{figure}

The situation is more interesting and, at the same time, more complicated, for connected graphs on $m=6$ vertices with a unique $1$-factor. To this end, we recall the well-known Kotzig's theorem stating that a graph with a unique 1-factor has a bridge that belongs to the 1-factor sub-graph. Splitting of $6$ vertices into two subsets of 3 vertices connected by a bridge leads to graphs $H_1, H_4, H_{19}$ shown in Fig.~\ref{fig-hexafamily}. Splitting into subsets of 2 and 4 vertices is impossible because the bridge should belong to the 1-factor and so the hanging leaf vertex of a 2-vertices sub-graph is not contained in the 1-factor. Splitting into a 1 vertex graph and 5-vertices graph lead to the remaining 17 graphs shown in Fig.~\ref{fig-hexafamily}. One can construct these 17 graphs from the set of all 10 graphs on four vertices (including disconnected graphs) by bridging to $K_2$ using up to $4$ edges. 

In summary, there exist 20 undirected connected graphs on $m=6$ vertices with a unique 1-factor shown in Fig.~\ref{fig-hexafamily}. All of them have invertible adjacency matrix. Except of the graph $H_{19}$ they are integrally invertible.

In this census, there are three bipartite graphs $H_1, H_2, H_6$ which are simultaneously positively and negatively invertible. There are twelve graphs 
\[
H_3, H_4, H_7, H_8, H_9, H_{13}, H_{14}, H_{15}, H_{16}, H_{17}, H_{18}, H_{20},
\] 
which are positively invertible. The three graphs $H_5, H_{10}, H_{12}$ are negatively invertible. The integrally invertible graph  $H_{11}$ is neither positively nor negatively invertible. The graphs $H_8$ and $H_{18}$ are iso-spectral but not isomorphic.

\begin{table}
\small
\caption{The family of graphs on 6 vertices with a unique 1-factor, their signability and spectrum. Graphs $H_8$ and $H_{18}$ are iso-spectral but not isomorphic.} 
\begin{center}
\begin{tabular}{c||c|c}
Graph & invertibility & spectrum \\
\hline\hline
   $H_1$ &  pos, neg & $\{ -1.8019, -1.2470, -0.4450, 0.4450, 1.2470, 1.8019 \}$ \\
\hline
   $H_2$ &  pos, neg & $\{ -1.9319, -1.0000, -0.5176, 0.5176, 1.0000, 1.9319 \}$ \\
\hline
   $H_3$ &  pos & $\{ -1.7397, -1.3738, -0.5945, 0.2742, 1.0996, 2.3342 \}$ \\
\hline
   $H_4$ &  pos & $\{ -1.7746, -1.0000, -1.0000, 0.1859, 1.3604, 2.2283 \}$ \\
\hline
   $H_5$ &  neg & $\{ -1.6180, -1.6180, -0.4142, 0.6180, 0.6180, 2.4142 \}$ \\
\hline
   $H_6$ &  pos, neg & $\{ -2.2470, -0.8019, -0.5550, 0.5550, 0.8019, 2.2470 \}$ \\
\hline
   $H_7$ &  pos & $\{ -1.8942, -1.3293, -0.6093, 0.3064, 0.7727, 2.7537 \}$ \\
\hline
   $H_8$ &  pos & $\{ -1.9032, -1.0000, -1.0000, 0.1939, 1.0000, 2.7093 \}$ \\
\hline
   $H_9$ &  pos & $\{ -1.6180, -1.3914, -1.0000, 0.2271, 0.6180, 3.1642 \}$ \\
\hline
$H_{10}$ &  neg & $\{ -1.8608, -1.6180, -0.2541, 0.6180, 1.0000, 2.1149 \}$ \\
\hline
$H_{11}$ &  int inv & $\{ -1.8241, -1.6180, -0.5482, 0.3285, 0.6180, 3.0437 \}$ \\
\hline
$H_{12}$ &  neg & $\{ -2.1420, -1.3053, -0.3848, 0.4669, 0.7661, 2.5991 \}$ \\
\hline
$H_{13}$ &  pos & $\{ -1.8563, -1.4780, -0.7248, 0.1967, 0.8481, 3.0143 \}$ \\
\hline
$H_{14}$ &  pos & $\{ -1.9202, -1.0000, -0.7510, 0.2914, 1.0000, 2.3799 \}$ \\
\hline
$H_{15}$ &  pos & $\{ -1.6783, -1.3198, -1.0000, 0.1397, 1.2297, 2.6287 \}$ \\
\hline
$H_{16}$ &  pos & $\{ -2.1364, -1.2061, -0.5406, 0.2611, 1.0825, 2.5395 \}$ \\
\hline
$H_{17}$ &  pos & $\{ -1.8619, -1.2827, -1.0000, 0.2512, 0.4897, 3.4037 \}$ \\
\hline
$H_{18}$ &  pos & $\{ -1.9032, -1.0000, -1.0000, 0.1939, 1.0000, 2.7093 \}$ \\
\hline
$H_{19}$ &  nonint inv & $\{ -1.7321, -1.0000, -1.0000, -0.4142, 1.7321, 2.4142 \}$ \\
\hline
$H_{20}$ &  pos & $\{ -2.3117, -1.0000, -0.6570, 0.3088, 0.7272, 2.9327 \}$ \\
\hline
\hline\hline
\end{tabular}

\end{center}
\noindent Legend: 'pos'/'neg' stands for a positively/negatively invertible graph, 'int inv' means an integrally invertible graph which is neither positively nor negatively invertible, 'nonint inv' stands for a graph with an adjacency matrix which is invertible but it is not integral.
\end{table}

\begin{table}
\label{tab-2}
\small
\caption{The family of graphs on 6 vertices with a unique 1-factor which can be arbitrarily bridged through $k=1,2,3$ vertices.} 
\begin{center}
\rotatebox{90}{
\begin{tabular}{c||c|c|c}
\tiny Graph &  $k=1$ & $k=2$ & $k=3$ \\
\hline\hline
\tiny   $H_1$ & \tiny $\{6\}, \{5\}, \{4\}, \{3\}, \{2\}, \{1\}$     &  \tiny $\{4,6\}, \{2,6\}, \{4,5\}, \{3,5\}, \{2,5\}, \{1,5\}, \{2,4\}, \{2,3\}, \{1,3\}$  & \tiny $\{2,4,6\}, \{2,4,5\}, \{2,3,5\}, \{1,3,5\}$  \\
\hline
\tiny   $H_2$ & \tiny $\{6\}, \{5\}, \{4\}, \{3\}, \{2\}, \{1\}$                 &  \tiny $\{4,6\}, \{2,6\}, \{3,5\}, \{2,5\}, \{1,5\}, \{3,4\}, \{2,4\}, \{1,4\}, \{2,3\}, \{1,3\}$  & \tiny $\{2,4,6\}, \{2,3,5\}, \{1,3,5\}, \{2,3,4\}, \{1,3,4\}$  \\
\hline
\tiny   $H_3$ & \tiny $\{6\}, \{5\}, \{4\}, \{3\}, \{2\}$ &  \tiny $\{4,6\}, \{2,6\}, \{3,5\}, \{2,5\}, \{3,4\}, \{2,4\}, \{2,3\}$   & \tiny $\{2,4,6\}, \{2,3,5\}, \{2,3,4\}$  \\
\hline
\tiny   $H_4$ & \tiny $\{6\}, \{5\}, \{4\}, \{2\}$ &  \tiny $\{4,6\}, \{2,6\}, \{4,5\}, \{2,5\}, \{2,4\}$  & \tiny $\{2,4,6\}, \{2,4,5\}$  \\
\hline
\tiny   $H_5$ & \tiny $\{6\}, \{5\}, \{4\}, \{3\}, \{2\}, \{1\}$            &  \tiny $\{3,6\}, \{2,6\}, \{4,5\}, \{3,5\}, \{2,5\}, \{1,5\}, \{2,4\}, \{2,3\}, \{1,3\}$   & \tiny $\{2,3,6\}, \{2,4,5\}, \{2,3,5\}, \{1,3,5\}$  \\
\hline
\tiny   $H_6$ & \tiny $\{6\}, \{5\}, \{4\}, \{3\}, \{2\}, \{1\}$               &  \tiny $\{5,6\}, \{4,6\}, \{2,6\}, \{3,5\}, \{1,5\}, \{3,4\}, \{2,4\}, \{2,3\}, \{1,3\}, \{1,2\}$          & \tiny $\{2,4,6\}, \{1,3,5\}, \{2,3,4\}, \{1,2,3\},$  \\
\hline
\tiny   $H_7$ & \tiny $\{5\}, \{4\}, \{3\}, \{2\}, \{1\}$             &  \tiny $\{3,5\}, \{1,5\}, \{3,4\}, \{2,4\}, \{2,3\}, \{1,3\}, \{1,2\}$            & \tiny $\{1,3,5\}, \{2,3,4\}, \{1,2,3\}$  \\
\hline
\tiny   $H_8$ & \tiny $\{4\}, \{3\}, \{2\}, \{1\}$            &  \tiny $\{3,4\}, \{2,4\}, \{1,4\}, \{2,3\}, \{1,2\}$            & \tiny $\{2,3,4\}, \{1,2,4\}$  \\
\hline
\tiny   $H_9$ & \tiny $\{4\}, \{3\}, \{2\}, \{1\}$               &  \tiny $\{3,4\}, \{2,4\}, \{2,3\}, \{1,3\}, \{1,2\}$        & \tiny $\{2,3,4\}, \{1,2,3\}$  \\
\hline
\tiny $H_{10}$ & \tiny $\{5\}, \{4\}, \{3\}, \{2\}, \{1\}$               &  \tiny $ \{4,5\}, \{2,5\}, \{3,4\}, \{2,4\}, \{1,4\}, \{1,3\}, \{1,2\}$       & \tiny $\{2,4,5\}, \{1,3,4\}, \{1,2,4\}$  \\
\hline
\tiny $H_{11}$ & \tiny $\{5\}, \{4\}, \{3\}, \{2\}, \{1\}$            &  \tiny $ \{4,5\}, \{2,5\}, \{3,4\}, \{2,4\}, \{1,4\}, \{1,3\}, \{1,2\}$   & \tiny $\{2,4,5\}, \{1,3,4\}, \{1,2,4\}$ \\
\hline
\tiny $H_{12}$ & \tiny $\{6\}, \{5\}, \{4\}, \{3\}, \{2\}, \{1\}$            &  \tiny $ \{5,6\}, \{4,5\}, \{2,5\}, \{3,4\}, \{2,4\}, \{1,4\}, \{1,3\}, \{1,2\}$          & \tiny $\{2,4,5\}, \{1,3,4\}, \{1,2,4\}$ \\
\hline
\tiny $H_{13}$ & \tiny $\{5\}, \{4\}, \{2\}, \{1\}$            &  \tiny $\{4,5\}, \{2,5\}, \{1,5\}, \{2,4\}, \{1,2\}$     & \tiny $\{2,4,5\}, \{1,2,5\}$\\
\hline
\tiny $H_{14}$ & \tiny $\{6\}, \{4\}, \{3\}, \{2\}, \{1\}$               &  \tiny $\{4,6\}, \{2,6\}, \{1,6\}, \{3,4\}, \{2,4\}, \{2,3\}, \{1,3\}, \{1,2\}$      & \tiny $\{2,4,6\}, \{1,2,6\}, \{2,3,4\}, \{1,2,3\}$ \\
\hline
\tiny $H_{15}$ & \tiny $ \{5\}, \{4\}, \{2\}, \{1\}$             &  \tiny $\{4,5\}, \{2,5\}, \{1,5\}, \{2,4\}, \{1,2\}$      & \tiny $\{2,4,5\}, \{1,2,5\}$ \\
\hline
\tiny $H_{16}$ & \tiny $\{6\}, \{5\}, \{3\}, \{2\}, \{1\}$             &  \tiny $\{5,6\}, \{3,5\}, \{1,5\}, \{2,3\}, \{1,3\}, \{1,2\}$        & \tiny $\{1,3,5\}, \{1,2,3\}$ \\
\hline
\tiny $H_{17}$ & \tiny $\{4\}, \{3\}, \{2\}, \{1\}$             &  \tiny $\{3,4\}, \{2,4\}, \{2,3\}, \{1,3\}, \{1,2\}$     & \tiny $\{2,3,4\}, \{1,2,3\}$\\
\hline
\tiny $H_{18}$ & \tiny $\{5\}, \{4\}, \{3\}, \{2\}, \{1\}$            &  \tiny $\{4,5\}, \{3,5\}, \{1,5\}, \{3,4\}, \{2,4\}, \{2,3\}, \{1,3\}, \{1,2\}$    & \tiny $ \{3,4,5\}, \{1,3,5\}, \{2,3,4\}, \{1,2,3\} $ \\
\hline
\tiny $H_{19}$ & ---                 &  ---              &  --- \\
\hline
\tiny $H_{20}$ & \tiny $\{6\}, \{4\}, \{3\}, \{2\}, \{1\}$                &  \tiny $\{4,6\}, \{1,6\}, \{3,4\}, \{2,4\}, \{2,3\}, \{1,3\}, \{1,2\}$     & \tiny $\{1,2,3\}, \{2,3,4\}$\\
\hline\hline
\end{tabular}
}
\end{center}

\end{table}

\section*{Acknowledgements} The authors are thankful to Prof.~Jozef~\v{S}ir\'a\v{n} and the anonymous referee for many helpful comments and suggestions that helped us to improve the paper. The research was supported by the Research Grants APVV 0136-12 and 15-0220, VEGA 1/0026/16 and 1/0142/17 (SP), and by the VEGA Research Grant 1/0780/15 (D\v{S}).

\end{document}